\documentclass[10pt]{article}

\usepackage{amsthm}
\usepackage{amsmath}
\usepackage{amsxtra}
\usepackage{amssymb}
\usepackage{thmtools}
\usepackage{enumitem}
\usepackage[retainorgcmds]{IEEEtrantools}
\usepackage{mathrsfs}
\usepackage{bbm}
\usepackage{accents}
\usepackage{graphicx}
\usepackage{turnstile}
\usepackage{tikz}
\usepackage[margin=1.2in]{geometry}
\usepackage{url}
\usepackage[colorlinks=true,linkcolor=blue]{hyperref}
\usepackage{comment}

\declaretheoremstyle[bodyfont=\sl]{slanted}

\swapnumbers
\declaretheorem[name=Definition,style=definition,qed=$\dashv$,
numberwithin=section]{dfn}
\declaretheorem[name=Theorem,style=slanted,sibling=dfn]{tm}

\declaretheorem[name=Remark,style=definition,sibling=dfn]{rem}

\declaretheorem[name=Fact,style=definition,sibling=dfn]{fact}

\swapnumbers
\declaretheoremstyle[headfont=\scshape]{claimstyle}
\declaretheorem[name=Claim,style=claimstyle]{clm}
\declaretheorem[name=Claim,style=claimstyle]{clm2}
\declaretheorem[name=Claim,style=claimstyle]{clm3}

\newcommand{\rest}{\upharpoonright}
\newcommand{\HOD}{\mathrm{HOD}}
\newcommand{\OD}{\mathrm{OD}}
\newcommand{\ZFC}{\mathrm{ZFC}}

\newcommand{\Coll}{\mathrm{Coll}}
\newcommand{\RR}{\mathbb{R}}
\newcommand{\OR}{\mathrm{OR}}
\newcommand{\pins}{\lhd}

\title{The cone fixing the  reals of $\HOD^{M_{2n}(x)}$}

\author{
  Farmer Schlutzenberg\footnote{Institute of Discrete Mathematics and Geometry, TU Vienna.\\\url{https://sites.google.com/site/schlutzenberg},
  email: \texttt{farmer dot schlutzenberg at tuwien dot ac dot at}}
  \and
  Connor Watson\footnote{\url{https://www.cjwsets.com}, email: \texttt{cjwsets@gmail.com}}
}

\begin{document}

\maketitle

\abstract{For a cone of reals $x$, the reals of $\mathrm{HOD}^{L[x]}$ are precisely the reals in $Q_3$, or equivalently, the set of reals in $M_1$, the minimal proper class mouse with a Woodin cardinal. We show that if $z$ is any real coding $M_1|\delta^{+M_1}$, then $z$ is the base of such a cone (so $M_1^{\#}$ is such a base). We also establish a kind of optimality for this result, by identifying a real $x$ such that $M_1\subseteq L[x]$ but $\mathbb{R}\cap\mathrm{HOD}^{L[x]}\not\subseteq M_1$, although there is a base for the cone in a generic extension of $L[x]$. We prove an analogous result for the reals of $\HOD^{M_{n}(x)}$ for even $n$, as well as an analogous optimality result which holds for all $n$.}

\section{The cone for $\RR\cap\mathrm{HOD}^{L[x]}$}

The set $Q_3$ of  descriptive set theory (more specifically, Q-theory) has many characterizations,
one being as the set of reals in $M_1$, the minimal proper class mouse with a Woodin cardinal.
In this note, we consider the following  characterization:

\begin{fact}\label{fact:cone}
For a Turing cone of reals $z$,
$Q_3=\RR\cap\HOD^{L[z]}$.
\end{fact}

The proof, given for example
in \cite[Proposition 2.41]{steel_cabal_part_5_intro},
involves the following two components:
\begin{enumerate}
\item\label{item:HOD^L[z]_sub} Genericity iterations are used to show
that for a Turing cone of reals $z$, every real in $\HOD^{L[z]}$
is in $M_1$,
\item\label{item:HOD^L[z]_sup} For all reals $z$ such that $M_1|\omega_1^{M_1}\subseteq L[z]$
and all $P\pins M_1|\omega_1^{M_1}$
such that $\rho_\omega^P=\omega$,
letting $\eta=\OR^{P}$, we have that
$L[z]\models$ ``$P$ is the unique $\omega$-sound, $(\omega,\omega_1+1)$-iterable premouse $P'$
such that $P'$ projects to $\omega$
and $\OR^{P'}=\eta$'',
and hence $P\in \HOD^{L[z]}$.
 \end{enumerate}

Since the fact above holds on a cone, 
it is natural to ask which reals might serve as a base for such a cone.
In the brief summary above, there is no specific base of the cone indicated in part \ref{item:HOD^L[z]_sub},
but note that $M_1^\#$ (or in fact
any real $z$ such that
$M_1|\omega_1^{M_1}\subseteq L[z]$) suffices 
as a base in part \ref{item:HOD^L[z]_sup}. So
one can hope
that $M_1^\#$ should serve as a base for both directions. We show here that this hope is fulfilled, and in fact:
\begin{tm}\label{tm:main}
Let $z\in\mathbb{R}$ be such that $M_1|\delta^{+M_1}\in L[z]$ and $\delta^{+M_1}$ is countable in $L[z]$, where $\delta$ is the Woodin cardinal of $M_1$. Then $\mathbb{R}\cap \HOD^{L[z]}=Q_3=\mathbb{R}\cap M_1$.
\end{tm}
So if $M_1^\#\leq_T z$ then the theorem applies to $z$.
Here ``$M_1^\#\leq_T z$'' literally means that $t\leq_T z$ where $t$ is the theory of $M_1^\#$.

We will actually give two proofs of the theorem,
as they provide somewhat different perspectives.

\begin{proof}
We first establish a claim which is
a more local version of the fact that for a cone of reals $x$, all reals
in $\HOD^{L[x]}$ are in $M_1$.

\begin{clm}\label{clm:cone_limiting_OD_reals}
There is a Turing cone of reals $x$ such
that whenever $\alpha$ is an ordinal and $y$ is a real in $L_\alpha[x]$ and
$L_\alpha[x] \models \ZFC^{-\varepsilon}$ + ``$y$ is OD'', then $y\in M_1$.\end{clm}

\begin{proof}If not then there is a Turing cone 
where it fails. For all such reals, let $\alpha_x$ be the
least  ordinal which is a counterexample;
that is, the least $\alpha$ such that $L_\alpha[x]\models \ZFC^{-\varepsilon}$ and there is a real $y \in\OD^{L_\alpha[x]}$ which is not in
$M_1$. Let $y_x$ be the resulting minimal $y$; note that $\alpha_x$ and
$y_x$ are Turing-invariant (though it is not immediate that they are constructibly invariant). Thus, for a Turing cone $C$ of reals $x$, we get the same $y_x$.

Now iterate $M_1$ to $N$ to make a base $b$ for $C$ generic for the
extender algebra. Let $g$ be $N[b]$-generic for $\Coll(\omega,\delta^N)$, so
$N[b,g]$ is equivalent to some $N[h]$ where h is a collapse generic, with
$b \leq_T x$ for  the natural real $x$ coding $(N|\delta^N,h)$. So $x\in C$,
so then $L_{\alpha_x}[x]$ defines $y_x$ from ordinal parameters, but then
by homogeneity of the
collapse, $y_x\in N$, so $y\in  M_1$, contradiction.
\end{proof}

Now working in $L[z]$ where $z$ is any real with $M_1|\delta^{+M_1}\in L[z]$ and countable in $L[z]$, let $g$ be
collapse generic over $M_1$ for collapsing $\omega_1^{M_1}$. Let $p=(M_1|\omega_1^{M_1},g)$. It is a
$\Sigma^1_3(\{p\})$ statement that ``there is a real $b$
such that for all reals $x\geq_T b$ and all ordinals $\alpha$ and all reals
$y \in L_\alpha[x]$, if $L_\alpha[x] \models \ZFC^{-\varepsilon}$ + ``$y$ is OD'', then $y\in 
M_1|\omega_1^{M_1}$. But then this statement has a witness in $M_1[g,h]$, whenever $h$ is $M_1[g]$-generic for collapsing $\delta$ to $\omega$ (seen by considering genericity iterations). Since $\delta^{+M_1}<\omega_1^{L[z]}$,
there is such an $h\in L[z]$,
and so we have such a base $b\in L[z]$.

Now let $y\in\RR\cap \HOD^{L[z]}$; we want to see that $y\in M_1$.
Let $z' = (b,z)$, so $z' \geq_T b$, and $L[z'] = L[z]$, so $y\in\HOD^{L[z']}$. But then there is $\alpha$ such that $L_\alpha[z']\models\ZFC^{-\varepsilon}$ + ``$y$ is OD'', and therefore $y \in M_1$.
\end{proof}
\begin{rem}
There is an alternate argument for part of the proof of Theorem \ref{tm:main},
which we now indicate. There is some real $b$ such that for all reals $y\geq_T b$,
all reals $w$ and all countable ordinals $\alpha$, if $w\in L_\alpha[y]\models\ZFC^{-\varepsilon}$ and $w\in \mathrm{OD}^{L_\alpha[y]}$ then $w\leq_T b$.
(This follows easily  from Claim \ref{clm:cone_limiting_OD_reals} above.) Moreover, this assertion is $\Sigma^1_3$.
So, in $M_1$, we have that $\mathrm{Coll}(\omega,\delta)$ forces that there is such a $b$. Take a pair of mutual generics $g_0$ and $g_1$ for $\mathrm{Coll}(\omega,\delta)$ which are computable from $M_1^\sharp.$ This gives witnesses $b_0$ and $b_1$ for the $\Sigma^1_3$ assertion just mentioned, with $b_i\in M_1[g_i]$
and $b_i\leq_T M_1^\sharp$. Now if
$y\geq_T M_1^\sharp$, then $y$ is also Turing above $b_0$ and $b_1.$ But then $y$ is computable from $b_0$ and $b_1$. So $y\in M_1[g_0]\cap M_1[g_1] = M_1.$
\end{rem}

\section{The cone for $\RR\cap \mathrm{HOD}^{M_{2n}(x)}$}
We would like to prove a generalization of Theorem \ref{tm:main}: if $x\geq_T M_n^\sharp$ then $\mathbb{R}\cap\mathrm{HOD}^{M_{n-1}(x)}=\mathbb{R}\cap M_n.$ Unfortunately, the methods used in the proof of Theorem \ref{tm:main} only work for odd $n$:

\begin{tm}\label{odd case of main tm}
Let $\delta_0$ be the least Woodin cardinal of $M_{2n+1}$. If $z$ is a real with $M_{2n+1}|\delta_0^{+M_{2n+1}}\in M_{2n}(z)$
and $\delta_0^{+M_{2n+1}}<\omega_1^{M_{2n}(z)}$ then $\mathbb{R}\cap \mathrm{HOD}^{M_{2n}(z)}  = \mathbb{R}\cap M_{2n+1}$.\end{tm}

\begin{proof}
    For one direction, suppose $a\in \mathbb{R}\cap M_{2n+1}$   . Choose an initial segment $P\trianglelefteq M_{2n+1}|\omega_1^{M_{2n+1}}$ with $a\in P$ and $\rho_\omega(P)=\omega$. Let $\eta=\mathrm{OR}^P$. Since $M_{2n+1}^\sharp\leq_T z$, we have $P\in M_{2n}(z).$ In $M_{2n}(z)$, $P$ is  characterized as the unique $(\omega,\omega_1+1)$-iterable  sound premouse $Q$ such that $\rho_\omega(Q)=\omega$ and $\mathrm{OR}^Q=\eta$.\footnote{One doesn't actually need to explicitly require full $(\omega,\omega_1+1)$-iterability here; we could replace this condition with $\Pi^1_{2n+2}$-iterability.} So $P\in\mathrm{HOD}^{M_{2n}(z)}$, and thus $a\in\mathrm{HOD}^{M_{2n}(z)}.$

    For the reverse inclusion, we prove a generalization of Claim \ref{clm:cone_limiting_OD_reals}:
 
    \begin{clm2}\label{clm:2n+1_Turing_cone}
        There is a Turing cone of reals $u$ such that whenever $\alpha$ is an ordinal and $y$ is a real with $y\in M_{2n}(u)|\alpha$ and $M_{2n}(u)|\alpha|\models\mathsf{ZFC}^{-\varepsilon}+``y\text{ is OD}",$ then $y\in M_{2n+1}$.\end{clm2}

        \begin{proof}
            Suppose not. Then for a Turing cone of reals $u$, there is a counterexample $\alpha$. For each such $u$, let $\alpha_u$ be the least ordinal $\alpha$ such that $M_{2n}(u)|\alpha\models\mathsf{ZFC}^{-\varepsilon}$
            and there is a real $y\in\mathrm{OD}^{M_{2n}(u)|\alpha}$ with $y\notin M_{2n+1}$. Let $y_u$ be the least real (in the $\OD$-order). There is then a cone $C$ and a real $y'$ such that $y_u=y'$ for all $u\in C.$ Let $b$ be a base for $C$. Genericity iterate $M_{2n+1}$ to an iterate $N$ such that $b$ is generic over $N$ for the extender algebra at $\delta_0^N$, the first Woodin cardinal of $N$.

            Let $g\subseteq\mathrm{Coll}(\omega,\delta_0^N)$ be $N[b]$-generic. The two-step extension $N[b,g]$ is equivalent to an extension $N[h]$ where $h$ is $N$-generic for $\mathrm{Coll}(\omega,\delta_0^N)$. Let $u$ be a real coding $(N|\delta_0^N,h)$.
            Then the extender sequence of $M_{2n}(u)$ is just given by the level-by-level extension of the sequence of $N$ above $\delta_0^N$, starting from base real $u$.
            So $M_{2n}(u)|\alpha_u$
            is just given by that process applied to $N|\alpha_u$.
            But $y'$ is definable over $M_{2n}(u)|\alpha_u$ from ordinal parameters.
            This shows that $y'$ is definable over the universe of $N[u]$ from parameters in $N$, and so 
            by the homogeneity of the collapse, $y'\in N$. But  $\mathbb{R}\cap N=\mathbb{R}\cap M_{2n+1}$, so $y'\in M_{2n+1}$, a contradiction.
        \end{proof}

Now fix a real $z$  with $M_{2n+1}|\delta_0^{+M_{2n+1}}\in M_{2n}(z)$
and $\delta_0^{+M_{2n+1}}<\omega_1^{M_{2n}(z)}$.
    We will show there is a base $b$ to the cone given by Claim  \ref{clm:2n+1_Turing_cone} with $b\in M_{2n}(z)$.
 Working in $M_{2n}(z)$, let $g$ be generic over $M_{2n+1}$ for collapsing $\omega_1^{M_{2n+1}}$, and set $p=(M_{2n+1}|\omega_1^{M_{2n+1}},g).$ It is a $\Sigma^1_{2n+3}(\{p\})$ assertion to say that there exists a real $b$ such that for all $u\geq_T b$, every ordinal $\alpha$, and all reals $y\in M_{2n}(u)|\alpha$, if $M_{2n}(u)|\alpha\models\mathsf{ZFC}^{-\varepsilon}+``y\text{ is OD}$'' then $y\in M_{2n+1}|\omega_1^{M_{2n+1}}$. (The key here is that we need only directly assert that for every $\Pi^1_{2n-1}$-iterable sound $(2n-2)$-small $u$-premouse $Q$ with $\rho_\omega^Q=\omega$, and every $\alpha\leq\mathrm{OR}^Q$, if $Q|\alpha\models\mathrm{ZFC}^{-\varepsilon}$ and $y\in\mathrm{OD}^Q$, then $y\in M_{2n+1}|\omega_1^{M_{2n+1}}$.) But then there is such a real  $b$ in $M_{2n+1}[g,h]$ whenever $h$ is $M_{2n+1}[g]$-generic for collapsing $\delta_0$ to $\omega$. 
 Therefore we get such a real $b\in M_{2n}(z)$, as desired.

  Now let $y$ be a real in $\OD^{M_{2n}(z)}$.
    We want to see that $y\in M_{2n+1}$. We can find an ordinal $\alpha>\omega_1^{M_{2n}(z)}$ such that $M_{2n}(z)\models\mathrm{ZFC}^{-\varepsilon}$ and $y\in \mathrm{OD}^{M_{2n}(z)|\alpha}$. 
     Let $\gamma<\omega_1^{M_{2n}(z)}$ with $b\in M_{2n}(z)|\gamma$ and $\rho_\omega^{M_{2n}(z)|\gamma}=
    \omega$. Let $u$ be the theory of $M_{2n}(z)|\gamma$. The extender sequence of $M_{2n}(u)$ is equivalent to that of $M_{2n}(z)$ above $\gamma$.
    But then $M_{2n}(u)|\alpha$ and $M_{2n}(z)|\alpha$ have the same universe, so $y\in\mathrm{OD}^{M_{2n}(u)|\alpha}$, and since $u\geq_T b$, therefore $y\in M_{2n+1}$, as desired.
\end{proof}

On the even side, fix $z\geq_T M_{2n+2}^{\#}$. Then we do have $\mathbb{R}^{M_{2n+2}}\subseteq \mathrm{OD}^{M_{2n+1}(z)}$,
for the same reasons as in the first paragraph of the foregoing proof. And Claim \ref{clm:2n+1_Turing_cone} also adapts directly. But the descriptive complexity and correctness calculations that follow that claim break down in the even context.

One could also consider ``$\mathrm{OD}^{L_{\omega_1^{\mathrm{ck}(z)}}(z)}$ for $z\geq_T 0^\#$."
In this case, we are not dealing with models of ZFC, so one first needs to determine what ``$\mathrm{OD}$'' should mean.
If one says that a real $y\in L_{\omega_1^{\mathrm{ck}(z)}}(z)$
is $\mathrm{OD}^{L_{\omega_1^{\mathrm{ck}(z)}}(z)}$ if there is some $\beta<\omega_1^{\mathrm{ck}(z)}$
and formula $\varphi$
such that for all reals $w\in L_{\omega_1^{\mathrm{ck}(z)}}(z)$, \[ y(m)=n \iff L_\beta[w]\models\varphi(w,m,n),\]
then straightforward absoluteness arguments show that any such $y$ is in $L$. But if one 
 fixes $0^\#$ at the base of the universe and says that $y$ is OD if there is $\beta<\omega_1^{\mathrm{ck}(z)}(z)$
such that
\[ y(m)=n\iff L_\beta\models\varphi(m,n),\]
then $0^\#$ itself is OD, and hence every real in $L_{\omega_1^{\mathrm{ck}(z)}}(z)$ becomes OD, so the notion becomes uninteresting. (Note that the formula $\varphi$ is in the language of set theory, so this fact is non-trivial; one cannot refer directly to the parameter $0^\#$. The point is that there are ordinals $\beta<L_{\omega_1^{\mathrm{ck}(z)}}(z)$ such that $0^\#$ is the unique iterable active sound premouse in $L_\beta[0^\#]$. For there is no other premouse which is fully iterable in $L_{\omega_1^{\mathrm{ck}(z)}}(z)$, since otherwise one easily gets $\omega_1^{\mathrm{ck}}(z)$-many (in ordertype) common indiscernibles,
and hence they generate the same theory of indiscernibles and are identical. And so a simple closure argument yields many ordinals $\beta$ in which every potential $0^\#$ in $L_\beta[z]$ has already been seen to be non-iterable.)

\section{Optimality of the hypotheses}
We now turn our attention to a kind of optimality result for Theorem \ref{odd case of main tm}.
 The main construction (giving the first paragraph of the theorem below) does not depend on the parity of the number of Woodin cardinals:

\begin{tm}\label{optimality result}
For every $n\in\omega$, there is a model of the form $M_n(x)$ with $x\in\mathbb{R}$,
such that $M_{n+1}\subseteq M_n(x)$,
\[\omega_1^{M_{n+1}}=\omega_1^{M_n(x)}<\delta_0^{+M_{n+1}}<\aleph_2^{M_n(x)},\]
and $\mathbb{R}\cap\mathrm{OD}^{M_n(x)}\not\subseteq M_{n+1}$.

Therefore, if $n$ is even
and $G$ is $M_n(x)$-generic
for $\mathrm{Coll}(\omega,\omega_1^{M_n(x)})$, then $M_n(x)[G]$ has the universe of $M_n(z)$ for a real $z$ satisfying
the conditions of Theorem \ref{odd case of main tm}.\end{tm}

So in the case $n$ is even, $x$ is only a collapse generic away from forming a base of a cone stabilizing the reals in $\mathrm{HOD}^{M_n(z)}$. 
We do not know of an example of an $x\in\RR$ with $\omega_1^{M_{n+1}}<\omega_1^{M_n(x)}$
and  $\RR\cap \mathrm{OD}^{M_n(x)}\not\subseteq M_{n+1}$.

\begin{proof}

    In \cite{JJ_new_con}, Jensen and Johnsbr\r{a}ten force over $L$ using Suslin trees to form a generic extension in which there is a non-constructible $\Delta^1_3$ real. It is straightforward to adapt this result to force over $M_{n+1}$ for any $n<\omega$.\footnote{One replaces levels $L_\alpha$ of $L$ with levels $M_{n+1}|\alpha$ of $M_{n+1}$ throughout. Otherwise, the only real place with a little difference is in \cite[Claim 6]{JJ_new_con}, where a basic condensation argument is carried through to show that the tree $T^*_{n+1}$ that they have constructed is Suslin in the model \emph{they} call $M_{n+1}$. Here one must observe that condensation for (what \emph{we} are calling) $M_{n+1}$ applies. But this is a standard consideration: The models $M$ of which we form hulls  are segments of $M_{n+1}|\omega_2^{M_{n+1}}$ with largest cardinal $\omega_1^M=\omega_1^{M_{n+1}}$, and so the transitive collapse maps $\pi:\bar{M}\to M$ have $\mathrm{crit}(\pi)=\omega_1^{\bar{M}}$ and $\pi(\omega_1^{\bar{M}})=\omega_1^M=\omega_1^{M_{n+1}}$. Condensation applies directly to each restriction $\pi\upharpoonright \bar{M}|\alpha$ where $\rho_\omega^{\bar{M}|\alpha}=\omega_1^{\bar{M}}$.} 
 This  gives a real $c$ such that $c$ is lightface definable in any proper class model $N$ with  $M_{n+1}[c]\subseteq N$ and $\omega_1^{M_{n+1}}=\omega_1^N$ and $N$ sufficiently closed under mice;
 for the latter, it is enough that $N$ can compute all the correct Q-structures to guide the $(\omega,\omega^N_1+1)$-strategies for $M_{n+1}|\omega_1^{M_{n+1}}$, as computed in $N$.
 We will be using $N$ of form $M_n(x)$,
 in which we will indeed have such closure under mice.
 
 It remains to force over $M_{n+1}[c]$ to produce a model of the form $M_n(x)$ for a real $x$, while keeping $c$ intact. First, force over $M_{n+1}[c]$ with $G\subseteq\mathrm{Coll}(\omega_1^{M_{n+1}},\delta_0^{+M_{n+1}})$,
 so that $\omega_1^{M_{n+1}[c][G]}<\delta_0^{+M_{n+1}}<\omega_2^{M_{n+1}[c][G]}$, and $\RR\cap M_{n+1}[c]=\RR\cap M_{n+1}[c][G]$. Then $M_{n+1}[c][G]$ has universe that of $M_n(X)$ where $X=(M_{n+1}|\omega_1^{M_{n+1}},c,G')$, represented as a subset of $\omega_1^{M_{n+1}}$, and where $G'$ is the natural coding of $(M_{n+1}|\delta_0^{+M_{n+1}}[c],G)$ as a subset of $\omega_1^{M_{n+1}}$.
 (The extender sequence of $M_n(X)$ is just that of $M_{n+1}$ above $\delta_0^{+M_{n+1}}$, extended to the generic extension.) Note that $\mathrm{HC}^{M_{n+1}[c]}=\mathrm{HC}^{M_{n+1}[c][G]}$ and every element of this set is encoded simply into $X$.
 
 In order to extend $M_{n+1}[c][G]$ into a model of the form $M_n(x)$ for an appropriate real $x$, we use a generalization of reshaping forcing, followed by almost disjoint forcing. 
Reshaping forcing is $\omega$-distributive, and   almost disjoint forcing is ccc, so the two-step extension preserves $\omega_1^{M_n(X)}$.

The reshaping forcing has been developed for mice at some places in the literature (for example, \cite{coding_into_K_reasonable_forcing_schindler}, 
 \cite{caicedo_schindler_projective_well_orders}), but we did not find its development at the level of $M_n$,
 so we discuss the details here.

    Let $\kappa = \omega_1^{M_{n+1}}=\omega_1^{M_n(X)}$. Work in $M_n(X)$.
    Let $\mathbb{P}$ be the forcing notion consisting of conditions $p$ such that for some limit $\gamma<\kappa$:
    \begin{enumerate}
        \item $p: \gamma\to 2$,
        \item for every limit $\lambda$ and every $k<\omega$ with $\lambda+2k<\gamma$, we have $p(\lambda+2k)=1$ iff $\lambda+k\in X$, and
        \item\label{item:reshaping_clause_3} for every nonzero limit $\xi\leq\gamma$, we have $M_n(p\upharpoonright\xi)\models$ ``$\xi$ is countable'',
    \end{enumerate}
    ordered by extension.
    Note that $\mathbb{P}$ is definable
    over the structure $(\mathrm{HC}^{M_{n}(X)},X)$,
    and therefore is in $M_n(X)|(\kappa+\omega)$,  since the mouse operator needed to specify clause \ref{item:reshaping_clause_3} is suitably definable over $\mathrm{HC}^{M_{n}(X)}$, since this set is closed under enough mice (that is, for each real $y\in M_n(X)$, we have $\mathbb{R}\cap M_n(y)\subseteq M_{n}(X)$).
    
    \begin{clm3}\label{clm:reshaping_omega-dist}
        $\mathbb{P}$ is $\omega$-distributive.
            \end{clm3}
        \begin{proof}
                    Fix $p\in \mathbb{P}$ and a sequence $\vec{D}=\langle D_k: k<\omega\rangle$ of dense open subsets $D_k\subseteq\mathbb{P}$.
                    Let $Q=M_n(X)|\omega_2^{M_n(X)}$. Let $Z=\mathrm{Hull}^Q(\{p,\vec{D}\})$ and $\alpha=Z\cap\kappa$. Fix $\langle\alpha_k:k<\omega\rangle$ with $\alpha=\sup_{k<\omega}\alpha_k$, and using elementarity, choose a sequence $\langle p_k:k<\omega\rangle$ recursively by setting $p_0=p$ and $p_{k+1}\in D_k\cap Z$ such that $p_{k+1}\leq p_k$ and $\mathrm{dom}(p_{k+1})>\alpha_k$. Then $q=\bigcup_{k<\omega} p_k$ is a condition beneath $p$ meeting all the $D_k$'s.
                    This is all quite immediate except for clause \ref{item:reshaping_clause_3} when $\xi=\alpha$; for this we have to see that $\alpha$ is countable in $M_n(q)$.
                    
                    So let $\bar{Z}$ be the transitive collapse of $Z$ and $\pi:\bar{Z}\to Z\preccurlyeq Q$ the uncollapse map. Then $\bar{Z}\pins M_n(X\cap\alpha)$, since $\bar{Z}$ is embedded into $Q$ (and hence $\bar{Z}$ is a stack of sound $n$-small $(X\cap\alpha)$-mice projecting to $\alpha$).
                    Therefore $\bar{Z}\in M_n(X\cap\alpha)$. But $\bar{Z}=\mathrm{Hull}^{\bar{Z}}(\{\pi^{-1}(p),\pi^{-1}(\vec{D})\})$, which determines a surjection of $\omega$ onto $\alpha$, so $\alpha$ is countable in $M_n(X\cap\alpha)$. So to see that $\alpha$ is countable in $M_n(q)$,
                    it suffices to see that $\mathcal{P}(\alpha)\cap M_n(X\cap\alpha)\subseteq M_n(q)$. 
                    But given any set $Y\subseteq\alpha$ with $Y\in M_n(X\cap\alpha)$, we can fix a sound $n$-small $(X\cap\alpha)$-mouse $P$ which projects $\leq\alpha$,  with $Y\in P$. Working in $M_n(q)$, 
                    since $\alpha<\delta_0^{M_n(q)}$, the maximal fully backgrounded $L[\mathbb{E},X\cap\alpha]$-construction produces a proper class model $N_\infty$ in which $\delta_i^{M_n(q)}$ is Woodin for each $i<n$
                    (where $\delta_i^{M_n(q)}$ is the $i$th Woodin of $M_n(q)$). Therefore $P\pins N_\infty$, and in particular, $P\in M_n(q)$, which suffices.
        \end{proof}

Let $g$ be $M_n(X)$-generic for $\mathbb{P}$. 
Since $\mathrm{HC}^{M_n(X)}=\mathrm{HC}^{M_{n+1}[c]}$ is coded simply into $X$,  $g$ is equivalent to a subset of $\kappa=\omega_1^{M_n(X)}$, modulo $X$.
Now $M_n(X)[g]$ is a forcing extension
of (the universe of) $M_n(X)$, so not formally a premouse, 
whereas $M_n(X,g)$ is formally an $(X,g)$-premouse.

\begin{clm3}\label{clm:M_n(X,g),M_n(X)[g]_have_same_univ} $M_n(X,g)$ and $M_n(X)[g]$ have the same universe.
\end{clm3}
\begin{proof}
This follows from the facts that $\kappa$ is a strong cutpoint of $M_n(X)$ and, as noted earlier, $\mathbb{P}$ is definable over $M_n(X)|\kappa\approx(\kappa,X)$. Thus,  the extender sequence of $M_n(X)$ (which has domain above $\kappa=\omega_1^{M_n(X)}$)
extends canonically to the generic extension via $g$, as the forcing is small relative to all of the critical points of extenders on the sequence. The resulting extender sequence is just that of $M_n(X,g)$. 
\end{proof}

The final step is almost disjoint forcing.
Let $f=\bigcup g$. Let  $\mathscr{F}=\left<A_\alpha\right>_{\alpha<\kappa}$ be the almost disjoint family on $\omega$ defined as follows. Given $\mathscr{F}\rest\eta$ for a fairly closed limit $\eta<\kappa$, $A_\beta$ is defined for $\beta\in[\eta,\lambda)$ for the next fairly closed ordinal $\lambda>\eta$, by the usual diagonalization process, using the least enumeration of $\eta$ in $M_n(f\rest\eta)$ to continue at stage $\eta$.
This enumeration exists by clause \ref{item:reshaping_clause_3} of the definition of $\mathbb{P}$.
Force with almost disjoint forcing with respect to the family $\mathscr{F}$, coding $f$ into the generic real. Let $h$ be the $M_n(X,g)$-generic filter and $x=x_h$ the generic real.

\begin{clm3}$M_n(X,g)[h]=M_n(X,g)[x]$ has universe that of $M_n(x)$.\end{clm3}
\begin{proof}Working in $M_n(x)$, we can recursively compute $\mathscr{F}\rest\eta$ and $f\rest\eta$ for $\eta\leq\kappa$, using $x$ and closure under the operator $y\mapsto\mathbb{R}\cap M_n(y)$, for reals $y\in M_n(x)$. So we get $\mathscr{F},f\in M_n(x)$. But $f$ encodes $X,g$, so $X,g\in M_n(x)$.

Now 
$\mathcal{P}(\kappa)\cap M_n(X,g)[h]=\mathcal{P}(\kappa)\cap M_n(x)$, like in the proof of Claim \ref{clm:reshaping_omega-dist}. So $\kappa=\omega_1^{M_n(X,g)[h]}=\omega_1^{M_n(x)}$ and $\omega_2^{M_n(X)}=\omega_2^{M_n(X,g)[h]}=\omega_2^{M_n(x)}$,
and $M_n(X,g)[h]|\omega_2^{M_n(X,g)[h]}$
has the same universe $R$ as does $M_n(x)|\omega_2^{M_n(x)}$.
Also, the extender sequence of $M_n(X,g)|\omega_2^{M_n(X,g)}$ is definable from parameters over $R$,
and the sequence of $M_n(x)|\omega_2^{M_n(x)}$ is also definable 
from parameters over $R$.
But since the almost disjoint forcing is in $M_n(X,g)|\omega_2^{M_n(X,g)}$
and $\omega_2^{M_n(X,g)}$ is a strong cutpoint
of $M_n(X,g)$,
the extender sequence of $M_n(X,g)$ above $\omega_2^{M_n(X,g)}$ extends canonically to the generic extension via $h$, giving a $M_n(X,g)[h]|\omega_2^{M_n(X,g)[h]}$-premouse $N$.
This therefore translates into an $M_n(x)|\omega_2^{M_n(x)}$-premouse $N'$.\footnote{So $N'$ has the same extender sequence as does $N$, which has critical points all above $\omega_2^{M_n(x)}$. The only difference between them is the structure at their bases. 
For $N$, that structure is $(R,\mathbb{E})$ for a certain predicate $\mathbb{E}$, and for $N'$, it is $(R,\mathbb{E}')$ for a certain predicate $\mathbb{E}'$.}
And now comparison of $M_n(x)$ versus $N'$ gives that $M_n(x)=N'$, and therefore $M_n(X,g)[h]$ has the same universe as does $M_n(x)$.
\end{proof}

So $M_{n+1}[c]\subseteq M_n(x)$, and since the forcings preserved $\omega_1^{M_{n+1}}=\omega_1^{M_{n+1}[c]}$, we have $\omega_1^{M_n(x)}=\omega_1^{M_{n+1}}$, which completes the proof.
\end{proof}

\section{Questions}

The following questions are the most important ones left open by the results here:
\begin{enumerate}
\item Is $M_2^\#$ the base of a cone
of reals $x$ such that $\RR\cap\HOD^{M_1(x)}\subseteq M_2$?
\item Is $\RR\cap\HOD^{M_1(M_2^\#)}\subseteq M_2$?
\item Is there a real $x$ such that $M_1\subseteq L[x]$ and $\omega_1^{M_1}<\omega_1^{L[x]}$ and $\mathbb{R}\cap\mathrm{HOD}^{L[x]}\not\subseteq M_1$?
\end{enumerate}

\section{Acknowledgements}

F.~Schlutzenberg supported by 
Austrian Science Fund (FWF) [PAT2418625].

Google Gemini was used for assistance with literature search in the preparation of this article.

\bibliographystyle{plain}
\bibliography{bibliography}
\end{document}